\documentclass[leqno]{article} 
\usepackage{latexsym}
\usepackage{amssymb}
\newtheorem{theorem}{Theorem}[section]

\newtheorem{remark}{Remark}[section]
\newtheorem{example}{Example}[section]

\newtheorem{proposition}{Proposition}[section]
\newtheorem{definition}{Definition}[section]
\begin{document}
\newcommand{\uno}{\hbox{\bf 1}}
\newcommand{\EE}{\mathop{\rm E}}
\newcommand{\remove}[1]{}
\newcommand{\barF}{\overline{F}}
\newcommand{\barH}{\overline{H}}
\newcommand{\finevuoto}{\hfill $\Box$}   
\newcommand{\qed}{\hfill\rule{2mm}{2mm}}   
\newenvironment{proof}{\begin{trivlist}
\item[\hspace{\labelsep}{\bf\noindent Proof. }]}
{\qed\end{trivlist}}
\title{\large\bf ON WEIGHTED RESIDUAL AND PAST ENTROPIES\footnote{
This is a corrected version of the paper appeared in [Sci.\ Math.\ Jpn.\ 64 (2006), no.\ 2, 255--266] and 
released on  June 21, 2011. The authors thank Fakhroddin Misagh for useful comments.}}
\author{\sc Antonio Di Crescenzo\\
\small\rm Dipartimento di Matematica e Informatica\\
\small\rm Universit\`a di Salerno\\
\small\rm Via Ponte don Melillo\\
\small\rm I-84084 Fisciano (SA), Italy\\
\small\tt adicrescenzo@unisa.it
\and
\sc Maria Longobardi\\
\small\rm Dipartimento di Matematica e Applicazioni\\
\small\rm Universit\`a di Napoli Federico II\\
\small\rm Via Cintia\\
\small\rm I-80126 Napoli, Italy\\
\small\tt maria.longobardi@unina.it
}

\date{\empty}

\maketitle
\begin{abstract}
We consider a ``length-biased'' shift-dependent information measure, 
related to the differential entropy in which higher weight is assigned to 
large values of observed random variables. This allows us to introduce the 
notions of ``weighted residual entropy'' and ``weighted past entropy'', 
that are suitable to describe dynamic information of random lifetimes, 
in analogy with the entropies of residual and past lifetimes introduced 
in \cite{Eb96} and \cite{DiLo2002}, respectively. The obtained results include 
their behaviors under monotonic transformations. 



\end{abstract}
%
\section{Introduction}
It is well known that knowledge and use of schemes for information coding and transmission 
play a relevant role in understanding and modeling certain aspects of biological systems 
features, such as neuronal activity. Since the pioneering contributions by Shannon 
\cite{Sh48} and Wiener \cite{Wi48} numerous efforts have been devoted to enrich and extend the 
underlying information theory. Various measures of uncertainty introduced in the 
past have been recently invoked in order to deal with information in the context of 
theoretical neurobiology (see, for instance, Johnson and Glantz \cite{JoGl04}). 
In addition, recent articles have thoroughly explored the use of information measures 
for absolutely continuous non-negative random variables, that appear to be suitable 
to describe random lifetimes (see \cite{EbPe95} and \cite{NaBeRuDe2004}, for instance). 
Here, wide use is made of Shannon entropy, that is also applied to residual 
and past lifetimes (cf., for instance, \cite{Eb96} and \cite{DiLo2002}). 
However, use of this type of entropy has the drawback of being position-free. 
In other terms, such an information measure does not take into account the values 
of the random variable but only its probability density. As a consequence, a random 
variable $X$ possesses the same Shannon entropy as $X+b$, for any $b\in\mathbb{R}$. 
\par
To come up with a mathematical tool whose properties are similar to those of 
Shannon entropy, however without being position-free, we introduce the 
notions of ``weighted residual  entropy'' and ``weighted past  entropy''. 
They are finalized to describe dynamic information of random lifetimes. 
In Section 2 we provide some basic notions on weighted entropy, complemented by some examples. 
Section 3 is devoted to the presentation of properties and of 
results on weighted residual and past  entropies, whose behaviors 
under monotonic transformations are studied in Section 4. 
\par
Throughout this paper, the terms ``increasing'' and ``decreasing'' are used in 
non-strict sense. Furthermore, we shall adopt the following notations: 
\begin{description}
\item{} $X$: an absolutely continuous non-negative honest random variable 
(representing for instance the random lifetime of a system or of a living organism, 
or the interspike intervals in a model of neuronal activity);
\item{} $f(x)$: the probability density function (pdf) of $X$;
\item{} ${\cal S}=(0,\nu),$ with $\nu\leq +\infty$: the support of $f(x)$; 
\item{} $F(x)=\Pr(X\leq x)$: the cumulative distribution function of $X$;
\item{} $\barF(x)=1-F(x)$: the survival function of $X$;
\item{} $\lambda(x)=f(x)/\barF(x)$: the hazard function, or failure rate, of $X$; 
\item{} $\tau(x)=f(x)/F(x)$: the reversed hazard rate function of $X$; 
\item{} $[Z|B]$: any random variable whose distribution is identical to the conditional 
distribution of $Z$ given $B$.
\end{description}
%
\section{Weighted entropy}
The differential entropy of $X$, or Shannon information measure, is defined by 
\begin{equation}
 H:=-\EE[\log f(X)]=-\int_0^{+\infty}f(x)\,\log f(x)\,{\rm d}x,
 \label{equation:1}
\end{equation}
where here ``$\log$'' means natural logarithm. The integrand function on the 
right-hand-side of (\ref{equation:1}) depends on $x$ only via $f(x)$, thus making 
$H$ shift-independent. Hence, $H$ stays unchanged if, for instance, $X$ is uniformly 
distributed in $(a,b)$ or $(a+h,b+h),$ whatever $h\in \mathbb{R}$. However, in 
certain applied contexts, such as reliability or mathematical neurobiology, it is 
desirable to deal with shift-dependent information measures. Indeed, knowing that 
a device fails to operate, or a neuron to release spikes in a given time-interval, 
yields a relevantly different information from the case when such an event occurs in 
a different equally wide interval. In some cases we are thus led to resort to a shift-dependent 
information measure that, for instance, assigns different measures to such distributions. 

In agreement with Belis and Guia\c{s}u \cite{BeGu1968} 
and Guia\c{s}u \cite{Gu86}, in this paper we shall refer to the following notion 
of {\em weighted entropy\/}:
\begin{equation}
 H^w:=-\EE[X\log f(X)]=-\int_0^{+\infty}x\,f(x)\,\log f(x)\,{\rm d}x,
 \label{equation:4}
\end{equation}
or equivalently:
$$
 H^w=-\int_0^{+\infty}{\rm d}y\int_y^{+\infty}f(x)\,\log f(x){\rm d}x.
$$
Recalling Taneja \cite{Ta90} we point out that the occurrence of an event removes a 
double uncertainty: a qualitative uncertainty related to its probability of occurrence, 
and a quantitative uncertainty concerning its value or its usefulness. The factor $x$, 
in the integral on the right-hand-side of (\ref{equation:4}), may be viewed as a 
weight linearly emphasizing the occurrence of the event $\{X=x\}$. 
This yields a ``length biased'' shift-dependent information measure assigning 
greater importance to larger values of $X$. The use of weighted entropy 
(\ref{equation:4}) is also motivated by the need, arising in various communication 
and transmission problems, of expressing the ``usefulness'' of events by means of an 
information measure given by $H^w=\EE[I(u(X),f(X))]$, where $I(u(x),f(x))$ satisfies
the following properties (see Belis and Guia\c{s}u \cite{BeGu1968}):
$$
 I(u(x),f(x)g(x))=I(u(x),f(x))+I(u(x),g(x)),
 \qquad 
 I(\lambda u(x),f(x))=\lambda I(u(x),f(x)).
$$
The relevance of weighted entropies as a measure of the average amount of 
{\em valuable\/} or {\em useful\/} information provided by a source 
has also been emphasized by Longo \cite{Lo76}.  
\par
The following are examples of pairs of distributions that possess same 
differential entropies but unequal weighted entropies. 
\begin{example}{\rm
Let $X$ and $Y$ be random variables with densities 
$$
 f_X(t)=\left\{
 \begin{array}{ll} 
 2t & \mbox{if $0<t<1$,} \\
 0 & \mbox{otherwise,} 
 \end{array}
 \right. \qquad
 f_Y(t)=\left\{
 \begin{array}{ll} 
 2(1-t) & \mbox{if $0< t< 1$,} \\
 0 & \mbox{otherwise.} 
 \end{array}
 \right.
$$
Their differential entropies are identical (see Example 1.1 
of \cite{DiLo2002}, where a misprint has to be noticed):
$$
 H_X= H_Y={1\over 2}-\log 2. 
$$ 
Hence, the expected uncertainties for $f_X$ and $f_Y$ on the 
predictability of the outcomes of $X$ and $Y$ are identical. Instead, $H_X^w < H_Y^w,$ since
$$
 H_X^w=-\int_0^1x\,2x\,\log 2x\,{\rm d}x=\frac{2}{9}-\frac{2}{3}\log 2,
$$
$$
 H_Y^w=-\int_0^1 2x\,(1-x)\log (1-x)\,{\rm d}x=\frac{1}{9}.
$$
Hence, even though $H_X=H_Y$, the expected
weighted uncertainty contained in $f_X$ on the predictability of the 
outcome of $X$ is smaller than that of $f_Y$ on the predictability of 
the outcome of $Y$. 
}\end{example}
\par
\begin{example}{\rm
Consider the piece-wise constant pdf 
$$
 f_X(x)=\sum_{k=1}^n c_k\,\uno_{\{k-1\leq x<k\}}
 \qquad \left(c_k\geq 0, \;\; k=1,2,\ldots,n; \quad \sum_{k=1}^n c_k=1\right).
$$ 
Its differential entropy is $H_X=-\sum_{k=1}^n c_k\log c_k$ while its weighted entropy is 
\begin{equation}
 H^w_X=-\sum_{k=1}^n k\,c_k\log c_k-{1\over 2}H_X.
 \label{equation:28}
\end{equation}
Note that any new density 
obtained by permutation of $c_1,c_2,\ldots,c_n$ has the same entropy 
$H_X$, whereas its weighted entropy is in general different from (\ref{equation:28}).
}\end{example}  
\par
Hereafter we recall some properties of differential entropy (\ref{equation:1}):\\
{\em (i)\/} for $a,b>0$ it is $H_{aX+b}=H_X+\log a$. \\
{\em (ii)\/} if $X$ and $Y$ are independent, then 
$$
 H_{(X,Y)}=H_X+H_Y,
$$ 
where   
$$
 H_{(X,Y)}:=-\EE[\log f(X,Y)] 
 =-\int_0^{+\infty}\int_0^{+\infty}f(x,y)\,\log f(x,y)\,{\rm d}x\,{\rm d}y
$$ 
is the bidimensional version of (\ref{equation:1}).
\par
Proposition 2.1, whose proof is omitted being straightforward, shows the corresponding properties
of $H^w$ similar to {\em (i)\/} and {\em (ii)\/}, but 
the essential differences emphasize the role of the mean value in the evaluation 
of the weighted entropy.
\begin{proposition}
The following statements hold:\\
{\em (i)\/} for $a,b>0$ it is $H_{aX+b}^w=a[H_X^w+\EE(X)\log a]+b(H_X+\log a)$. \\
{\em (ii)\/} if $X$ and $Y$ are independent, then  
$$
 H_{(X,Y)}^w=\EE(Y)\,H_X^w+\EE(X)\,H_Y^w,
$$ 
where $H_{(X,Y)}^w:=-\EE[X\,Y\,\log f(X,Y)]$. 
\end{proposition}
\par
Let us now evaluate the weighted entropy of some random variables.
\begin{example} 
{\rm {\it (a)} $X$ is exponentially distributed with parameter $\lambda>0.$ Then,
\begin{equation}
 H^w=-\int_0^{+\infty}x\,\lambda e^{-\lambda\,x}\,
 \log \left(\lambda e^{-\lambda\,x}\right) {\rm d}x 
 =\frac{2-\log \lambda}{\lambda}.
 \label{equation:20}
\end{equation}
{\it (b)} $X$ is uniformly distributed over $[a,b]$. Then
\begin{equation}
 H^w=-\int_a^b x\,\frac{1}{b-a}\,\log \frac{1}{b-a}\,{\rm d}x
 =\frac{a+b}{2}\,\log(b-a).
 \label{equation:21}
\end{equation}
It is interesting to note that in this case the weighted entropy can be expressed 
as the product 
\begin{equation}
 H^w=\EE(X)\,H.
 \label{equation:29}
\end{equation}
{\it (c)} $X$ is Gamma distributed with parameters $\alpha$ and $\beta$: 
$$
 f(x)=\left\{
 \begin{array}{ll} 
 \displaystyle\frac{x^{\alpha-1}e^{-x/\beta}}{\beta^{\alpha}\Gamma(\alpha)} 
 & \mbox{if $x>0$,}
 \\
 0 & \mbox{otherwise.} 
 \end{array}
 \right.
$$
There holds: 
\begin{eqnarray*}
 H^w \!\!\! & = & \!\!\! -\int_0^{+\infty} 
 \frac{x^{\alpha} e^{-\frac{x}{\beta}}}{\beta^{\alpha}\Gamma(\alpha)} \,
 \log \frac{x^{\alpha-1}e^{-\frac{x}{\beta}}}{\beta^{\alpha}\Gamma(\alpha)} \,{\rm d}x \\
 & = & \!\!\! \alpha \beta \,\log[\beta^{\alpha}\Gamma(\alpha)] 
 -\alpha(\alpha-1)\beta\,\big\{\log\beta+ \psi_0(\alpha+1)\big\}+\alpha(\alpha+1)\beta,
\end{eqnarray*}
where 
\begin{equation}
 \psi_0(x):=\frac{\rm d}{{\rm d}x}\log \Gamma(x)=\frac{\Gamma'(x)}{\Gamma(x)} 
 \label{equation:5}
\end{equation}
and where we have used the following identities, with $\alpha>-1$ and $\beta>0$:
$$
 \int_0^{+\infty} x^{\alpha}
 e^{-\frac{x}{\beta}}\,{\rm d}x=\left({\frac{1}{\beta}}\right)^{-1-\alpha}\;\Gamma(1+\alpha),
$$
$$
 \int_0^{+\infty} x^{\alpha}
 e^{-\frac{x}{\beta}}\,\log x\,{\rm d}x
 =\left({\frac{1}{\beta}}\right)^{-1-\alpha}\,
 \Gamma(1+\alpha)\,\big\{\log \beta+ \psi_0(1+\alpha)\big\}.
$$
{\it (d)} $X$ is Beta distributed with parameters $\alpha$ and $\beta$:
$$
 f(x)=\left\{
 \begin{array}{ll} 
 \displaystyle\frac{x^{\alpha-1} (1-x)^{\beta-1}}{B(\alpha,\beta)} 
 & \mbox{if $0<x<1$,}
 \\
 0 & \mbox{otherwise,} 
 \end{array}
 \right.$$
with 
$$
 B(\alpha,\beta)=\int_0^1x^{\alpha-1} (1-x)^{\beta-1}\,{\rm d}x
 =\frac{\Gamma(\alpha)\Gamma(\beta)}{\Gamma(\alpha+\beta)}.
$$
Then, 
\begin{eqnarray*}	
 H^w \!\!\!  
 &=& \!\!\! \frac{\log
 B(\alpha,\beta)}{B(\alpha,\beta)}\frac{\Gamma(\alpha+1)\Gamma(\beta)}
 {\Gamma(\alpha+\beta+1)}
 -\frac{(\alpha-1)}{B(\alpha,\beta)}\frac{\Gamma(\alpha+1)\Gamma(\beta)}
 {\Gamma(\alpha+\beta+1)}\,\big\{\psi_0(\alpha+1)-\psi_0(\alpha+\beta+1)\big\}  \\
 &-& \!\!\! 
 \frac{(\beta-1)}{B(\alpha,\beta)}\frac{\Gamma(\beta)}
 {\Gamma(\alpha+\beta+1)}\,\big\{\psi_0(\beta)\Gamma(\alpha+1)-\alpha
 \Gamma(\alpha)\psi_0(\alpha+\beta+1)\big\},
\end{eqnarray*}
where $\psi_0(x)$ is defined in (\ref{equation:5}). 
}\end{example}
\par
\begin{remark}{\rm 
We point out that different distributions may have identical weighted entropies. 
For instance, if $X$ is uniformly distributed over $[0,1]$ and $Y$ is exponentially 
distributed with mean $e^{-2}$, then from (\ref{equation:20}) and (\ref{equation:21}) 
we have 
$$
 H_X^w=H_Y^w=0.
$$
Another example is the following: if $X$ is uniformly distributed in $[0,2]$ 
and $Y$ is exponentially distributed with parameter $\lambda=1{.}93389$ 
(such that $\log \lambda+ \lambda \log 2=2$), then from (\ref{equation:20}) 
and (\ref{equation:21}) one obtains: 
$$
 H_X^w=H_Y^w=\log 2.
$$
}\end{remark}
\begin{remark}{\rm 
There exist random variables having negative arbitrarily large weighted entropy. 
For instance, if $X$ is uniformly distributed over $[a,b]$, with $b>0$, 
then from (\ref{equation:21}) we have   
$$
 \lim_{a\to b^-}H^w = \lim_{a\to b^-}\frac{a+b}{2}\,\log(b-a) 
 =-\infty. 
$$
}\end{remark}
\begin{remark}{\rm 
Notice that in general $H^w$ can be either larger or smaller than $H$. For instance, 
if $X$ is uniformly distributed over $[a,b]$, from (\ref{equation:29}) it 
follows $H^w\geq H$ when $\EE(X)\geq 1$, and $H^w\leq H$ when $\EE(X)\leq 1$. 
}\end{remark}
\begin{remark}{\rm 
If the support of $X$ is $[0,\nu]$, with $\nu$ finite, the following upper bound 
for the weighted entropy of $X$ holds:  
\begin{equation}
 H^w\leq \mu\log {\nu^2\over 2\mu}, 
 \qquad \hbox{where $\mu=\EE(X)\in[0,\nu]$.}
 \label{equation:24}
\end{equation}
This follows via the continuous version of Jensen's inequality. 
The maximum of $b(\mu):= \mu\log {\nu^2\over 2\mu}$ is attained at $\mu=\mu_M$, where
$$
 \mu_M=\left\{
 \begin{array}{ll}
 \displaystyle{\nu^2\over 2e} & \hbox{if $\nu<2e$,} \\
 \nu & \hbox{if $\nu\geq 2e$,} 
 \end{array}
 \right.
 \qquad \hbox{with}
 \quad 
 b(\mu_M)=\left\{
 \begin{array}{ll}
 \displaystyle{\nu^2\over 2e} & \hbox{if $\nu<2e$,} \\
 \nu\log\displaystyle{\nu\over 2} & \hbox{if $\nu\geq 2e$.} 
 \end{array}
 \right.
$$
If, in particular, $X$ is uniformly distributed over $[0,\nu]$, then (\ref{equation:24}) 
holds with the equal sign since $H^w=\frac{\nu}{2}\,\log\nu$ and $\mu=\nu/2$. 
}\end{remark}
%
\section{Weighted residual and past  entropies}
The residual entropy at time $t$ of a random lifetime $X$ was introduced by 
Ebrahimi \cite{Eb96} and defined as: 
\begin{eqnarray}	
 H(t) \!\!\! 
 &=& \!\!\! -\int_t^{+\infty}\frac{f(x)}{\barF(t)}\log \frac{f(x)}{\barF(t)}\,{\rm d}x 
 \nonumber \\
 &=& \!\!\! \log \barF(t)-\frac{1}{\barF(t)}\int_t^{+\infty}f(x)\log f(x)\,{\rm d}x  
 \nonumber \\
 &=& \!\!\! 1-\frac{1}{\barF(t)}\int_t^{+\infty}f(x)\log \lambda(x)\,{\rm d}x, 
 \label{equation:2}
\end{eqnarray}
for all $t\in {\cal S}$. We note that $H(t)$ is the differential entropy 
of the residual lifetime of $X$ at time $t$, i.e. of $[X\,|\,X>t]$. 
Various result on $H(t)$ have been the object of recent researches (see   
\cite{AsEb2000}, \cite{BeGuNaRu2001}, \cite{Eb97}, \cite{Eb2000}, \cite{EbKi96a}).  
The past entropy at time $t$ of $X$ is defined by Di Crescenzo and 
Longobardi \cite{DiLo2002} as:
\begin{equation}
 \barH(t)=-\int_0^t\frac{f(x)}{F(t)}\log\frac{f(x)}{F(t)}\,{\rm d}x\,, 
 \qquad t\in {\cal S},
 \label{equation:3}
\end{equation}
where $f(x)/F(t)$, $0<x<t$, is the pdf of the past lifetime $[X\,|\,X\leq t]$. 
Hence, given for instance that an item has been found failing at time $t$, $\barH(t)$ 
measures the uncertainty about its past life. We remark that $\barH(t)$ can 
also be viewed as the entropy of the  inactivity time $[t-X\,|\,X\leq t]$. 
\par
Various dynamic information functions have been recently introduced to measure 
discrepancies between residual lifetime distributions \cite{EbKi96b} and between 
past lifetime distributions \cite{DiLo04}, as well as to measure dependence between 
two residual lifetimes \cite{DiLoNa04}. In order to introduce new shift-dependent 
dynamic information measures, we now make use of (\ref{equation:4}) to define 
two weighted entropies for residual lifetimes and past lifetimes that are the 
weighted version of entropies (\ref{equation:2}) and (\ref{equation:3}). 
\begin{definition}
For all $t\in {\cal S}$, \\
(i) the weighted residual  entropy at time $t$ of a random lifetime $X$ 
is the differential weighted entropy of $[X\,|\,X>t]$:
\begin{equation}
 H^w(t):=-\int_t^{+\infty}x\,\frac{f(x)}{\barF(t)}\log \frac{f(x)}{\barF(t)}\,{\rm d}x; 
 \label{equation:6}
\end{equation}
(ii) the weighted past  entropy at time $t$ of a random lifetime $X$ 
is the differential weighted entropy of $[X\,|\,X\leq t]$:
\begin{equation}
 \barH^w(t):=-\int_0^t x\,\frac{f(x)}{F(t)}\log \frac{f(x)}{F(t)}\,{\rm d}x. 
 \label{equation:7} 
\end{equation}
\end{definition}
\par
We notice that 
$$
 \lim_{t\to +\infty}H^w(t)=\lim_{t\to 0^+}\barH^w(t)=H^w.
$$
In addition, due to (\ref{equation:6}), the weighted residual  entropy can be rewritten as 
\begin{equation}
 H^w(t)=-\frac{1}{\barF(t)}\int_t^{+\infty}x\,f(x)\log f(x)\,{\rm d}x
 +\frac{\log\barF(t)}{\barF(t)}\int_t^{+\infty}x\,f(x)\,{\rm d}x.
 \label{equation:8}
\end{equation}
The second integral in (\ref{equation:8}) can be calculated by noting that  
\begin{equation}
 \int_t^{+\infty}x\,f(x)\,{\rm d}x
 =t\,\barF(t)+\int_t^{+\infty}\barF(y)\,{\rm d}y.
 \label{equation:9}
\end{equation}
Furthermore, an alternative way of writing $H^w(t)$ is the following: 
\begin{eqnarray*}	
 H^w(t) \!\!\!  
 &=& \!\!\!  -\int_t^{+\infty}{\rm d}x\int_0^x \frac{f(x)}{\barF(t)}\log \frac{f(x)}
 {\barF(t)}\,{\rm d}y  \\
 &=& \!\!\!  -\int_0^t{\rm d}y\int_t^{+\infty} \frac{f(x)}{\barF(t)}\log
 \frac{f(x)} {\barF(t)}\,{\rm d}x-\int_t^{+\infty}{\rm d}y\int_y^{+\infty}
 \frac{f(x)}{\barF(t)}\log\frac{f(x)}{\barF(t)}\,{\rm d}x.  
\end{eqnarray*}
Recalling (\ref{equation:2}), we thus obtain: 
\begin{equation}
 H^w(t)=t\,H(t)
 +\frac{1}{\barF(t)}\int_t^{+\infty}\barF(y)\left[H(y)+ \log\frac{\barF(t)}{\barF(y)}\right]\,{\rm d}y.
 \label{equation:10}
\end{equation}
From (\ref{equation:10}) it easily follows
\begin{eqnarray}
 \frac{{\rm d}}{{\rm d}t}H^w(t) \!\!\!  
 &=& \!\!\! t\,\frac{{\rm d}}{{\rm d}t}H(t)+\frac{f(t)}{[{\barF(t)}]^2}\bigg[ \int_t^{+\infty}\barF(y)H(y)\,{\rm d}y 
 \nonumber \\
 & & \!\!\! - [1-\log \barF(t)]\int_t^{+\infty}\barF(y)\,{\rm d}y  -\int_t^{+\infty}\barF(y)\log \barF(y)\,{\rm d}y\bigg].
 \label{equation:11}
\end{eqnarray}
\par
In analogy with Theorem 1 of Belzunce {\em et al.}\ \cite{NaBeRuDe2004}, the
following characterization result holds. 
\begin{theorem}
If $X$ has an absolutely continuous distribution function $F(t)$ and if $H(t)$ is 
increasing for all $t\in {\cal S}$, then $H^w(t)$ uniquely determines $\barF(t)$. 
\end{theorem}
\begin{proof}
From (\ref{equation:8}) and (\ref{equation:9}) we have: 
$$
 \int_t^{+\infty}x\,f(x)\log f(x)\,{\rm d}x
 =\log\barF(t)\left[t\,\barF(t)+\int_t^{+\infty}\barF(y)\,{\rm d}y\right]
 -H^w(t)\,\barF(t).
$$
Differentiating both sides we obtain:
$$
 -t\,f(t)\,\log f(t)
 =-t\,f(t)\,\log \barF(t)-t\,f(t)-\frac{f(t)}{\barF(t)}\int_t^{+\infty}\barF(y)\,{\rm d}y
 +f(t)\,H^w(t)-\barF(t)\,\frac{\rm d}{{\rm d}t}H^w(t).
$$
Hence, due to (\ref{equation:11}) and recalling the hazard function 
$\lambda(t)=f(t)/\barF(t)$, it follows:
$$
 t \lambda(t)[1-\log \lambda(t)]  
 = \lambda(t)H^w(t)-t\,\frac{{\rm d}}{{\rm d}t}H(t) 
-\frac{\lambda(t)}{\barF(t)}\,{\cal I}(t),
$$
where
$$
 {\cal I}(t)=\int_t^{+\infty}\barF(y)H(y)\,{\rm d}y 
 + \log \barF(t)\int_t^{+\infty}\barF(y)\,{\rm d}y  -\int_t^{+\infty}\barF(y)\log \barF(y)\,{\rm d}y. 
$$
Then, for any fixed $t\in {\cal S}$, $\lambda(t)$ is a positive solution of equation $g(x)=0$, 
where 
\begin{equation}
 g(x) := x \left[t\,(1-\log x)- H^w(t)+\frac{1}{\barF(t)}\,{\cal I}(t)\right]
+t\,\frac{{\rm d}}{{\rm d}t}H(t).
 \label{equation:12}
\end{equation}
Note that $\displaystyle\lim_{x\to+\infty} g(x)=-\infty$ and 
$g(0)=t\,\frac{\rm d}{{\rm d}t}H(t)\geq 0$. Furthermore, from (\ref{equation:12}) we have:
$$
 \frac{\rm d}{{\rm d}x}g(x)
 =-t\,\log x-H^w(t)+ \frac{1}{\barF(t)}\,{\cal I}(t),
$$
so that $\frac{\rm d}{{\rm d}x}g(x)=0$ if and only if
$$
 x=\exp\left\{-\frac{1}{t}\left[H^w(t)
 -\frac{1}{\barF(t)}\,{\cal I}(t)\right]\right\}.
$$ 
Therefore, $g(x)=0$ has a unique positive solution so that $\lambda(t)$, and 
hence $\barF(t)$, is uniquely determined by $H^w(t)$ under assumption 
$\frac{\rm d}{{\rm d}t}H(t)\geq 0$. This concludes the proof. 
\end{proof}
\begin{remark}{\rm 
We note that, by virtue of (\ref{equation:11}), if $H(t)$ is increasing [decreasing], 
then also $H^w(t)$ is increasing [decreasing].
}\end{remark}
\par
In order to attain a decomposition of the weighted entropy, similar to that 
given in Proposition 2.1 of Di Crescenzo and Longobardi \cite{DiLo2002}, 
for a random lifetime $X$ possessing finite mean $\EE(X)$ we recall that  
the {\em length-biased distribution function\/} and the {\em length-biased 
survival function\/} are defined respectively as 
\begin{equation}
 F^*(t)=\frac{1}{\EE(X)}\int_0^t x\,f(x)\,{\rm d}x,
 \qquad \barF^*(t)=\frac{1}{\EE(X)}\int_t^{+\infty} x\,f(x)\,{\rm d}x,
 \qquad t\in{\cal S}.
 \label{equation:22}
\end{equation}
These functions characterize weighted distributions that arise in sampling 
procedures where the sampling probabilities are proportional to sample 
values. (See Section 3 of Belzunce {\em et al.}\ \cite{NaBeRuDe2004} for 
some results on uncertainty in length-biased distributions. Moreover, 
see Navarro {\em et al.}\ \cite{NaDeRu2001}, Bartoszewicz and 
Skolimowska \cite{BaSk2006} and references therein for characterizations 
involving weighted distributions). 
\begin{theorem}
For a random lifetime $X$ having finite mean $\EE(X)$, for all $t\in {\cal S}$ 
the weighted entropy can be expressed as follows:
$$
 H^w = \EE(X) \left\{-F^*(t)\,\log F(t)-\barF^*(t)\,\log\barF(t)\right\}
 +F(t)\,\barH^w(t)+\barF(t)\,H^w(t). 
$$
\end{theorem}
\begin{proof}
Recalling Eqs.\ (\ref{equation:4}), (\ref{equation:6}) and (\ref{equation:7}) we have:
\begin{eqnarray*}
 H^w \!\!\! &=& \!\!\! - F(t)\int_0^t x\,\frac{f(x)}{F(t)}\log f(x)\,{\rm d}x
 -\barF(t) \int_t^{+\infty}x\,\frac{f(x)}{\barF(t)}\log f(x)\,{\rm d}x \\
 &=& \!\!\! -\log F(t) \int_0^t x \,f(x)\,{\rm d}x-\log\barF(t)\int_t^{+\infty}x\,f(x)\,{\rm d}x
 +F(t)\,\barH^w(t)+\barF(t)\,H^w(t).
\end{eqnarray*}
The proof then follows from (\ref{equation:22}). 
\end{proof}
\par
In order to provide a lower bound for the weighted residual  entropy of a 
random lifetime $X$, let us introduce the following conditional mean value:
\begin{equation}
 \delta(t):=\EE(X\,|\,X>t)
 =\frac{1}{\barF(t)}\int_t^{+\infty}x\,f(x)\,{\rm d}x
 =t+\frac{1}{\barF(t)} \int_t^{+\infty} \barF(x)\,{\rm d}x,
 \qquad t\in{\cal S}.
 \label{equation:23} 
\end{equation}
\begin{theorem}
If the hazard function $\lambda(t)$ is decreasing in $t\in {\cal S}$, then
\begin{equation}
 H^w(t)\geq - \delta(t)\,\log \lambda(t), 
 \qquad t\in {\cal S}.
 \label{equation:14} 
\end{equation}
\end{theorem}
\begin{proof}
From (\ref{equation:6}) we have:
$$
 H^w(t)= -\frac{1}{\barF(t)}\int_t^{+\infty}x\,f(x)\,\log\lambda(x)\,{\rm d}x
 -\frac{1}{\barF(t)}\int_t^{+\infty}x\,f(x)\,\log\frac{\barF(x)} {\barF(t)}\,{\rm d}x, 
 \qquad t\in {\cal S}.
$$
Since $\log \frac{\barF(x)}{\barF(t)}\leq 0$ for $x\geq t$ and, 
by assumption, $\log \lambda(x) \leq \log \lambda(t)$, there holds:
\begin{eqnarray*}
 H^w(t) \!\!\! 
 &\geq & \!\!\! -\frac{1}{\barF(t)}\int_t^{+\infty}x\,f(x)\,\log\lambda(x)\,{\rm d}x \\
 &\geq & \!\!\! -\frac{\log \lambda(t)}{\barF(t)}\int_t^{+\infty}x\,f(x)\,{\rm d}x.
\end{eqnarray*}
The proof then follows by recalling (\ref{equation:23}).
\end{proof}
\par
In the following example we consider the case of constant hazard function. 
\begin{example}{\rm 
For an exponential distribution with parameter $\lambda>0$, the weighted residual 
entropy is given by
\begin{equation}
 H^w(t)=-\int_t^{+\infty} x\, \frac{\lambda e^{-\lambda\,x}}{e^{-\lambda\,t}}\,
 \log\frac{\lambda e^{-\lambda\,x}}{e^{-\lambda\,t}}\,{\rm d}x
 =t+\frac{2}{\lambda}-\left(t+\frac{1}{\lambda}\right)\log \lambda, 
 \qquad t\geq 0.
 \label{equation:27}
\end{equation}
Recalling that for an exponential r.v.\ $\lambda(t)=\lambda$ and 
$\delta(t)=t+\frac{1}{\lambda}$, it is easily seen that (\ref{equation:14}) is fulfilled. 
}
\end{example}
\par
Let us now discuss some properties of the weighted past  entropy.
From (\ref{equation:7}) we have: 
$$
 \barH^w(t)=-\frac{1}{F(t)}\int_0^t x\,f(x)\log f(x)\,{\rm d}x
 +\frac{\log F(t)}{F(t)}\int_0^tx\,f(x)\,{\rm d}x, 
 \qquad t\in {\cal S},
$$
where, similarly to (\ref{equation:9}), it  is  
\begin{equation}
 \int_0^t x\,f(x)\,{\rm d}x=t\,F(t)-\int_0^tF(y)\,{\rm d}y. 
 \label{equation:13}
\end{equation}
Alternatively,  
\begin{eqnarray}
 \barH^w(t)  \!\!\! 
 &=& \!\!\! -\int_0^t{\rm d}x\int_0^x \frac{f(x)}{F(t)}\log \frac{f(x)}{F(t)}\,{\rm d}y
 \nonumber \\
 &=& \!\!\! t\,\barH(t)-\frac{1}{F(t)}\int_0^t F(y)\left[\barH(y)+\log \frac{F(t)}{F(y)}\right]\,{\rm d}y,
 \label{equation:16}
\end{eqnarray}
where $\barH(t)$ is the past entropy given in (\ref{equation:3}).
\begin{example}{\rm 
The weighted past  entropy of an exponentially distributed random variable 
with parameter $\lambda>0$ is instead given by 
$$
 \barH^w(t)
 =\frac{1}{1-e^{-\lambda t}}
 \left[\frac{2}{\lambda} - \frac{2}{\lambda}e^{-\lambda t}
 - 2 t\,e^{-\lambda t}-\lambda t^2\,e^{-\lambda t}
 +\left(\frac{1}{\lambda} - \frac{1}{\lambda}e^{-\lambda t}-t\,e^{-\lambda t}\right)
 \log \frac{1-e^{-\lambda t}}{\lambda}\right] 
$$
for $t>0$. 
}\end{example}
\par
In order to obtain upper bounds for the weighted past  entropy let us now 
recall the definition of mean past lifetime:
\begin{equation}
 \mu(t):=\EE(X|X\leq t)
 =\int_0^t x \,\frac{f(x)}{F(t)}\,{\rm d}x=t-\frac{1}{F(t)}\int_0^t F(y)\,{\rm d}y.
 \label{equation:15} 
\end{equation}
We incidentally note that this is related to the Bonferroni Curve by   
$\mu(t)=\EE(X)\,B_F[F(t)]$ (see Giorgi and Crescenzi \cite{GiCr01}). 
\begin{theorem}
(i) For all $t\in {\cal S}$, it is
\begin{equation}
 \barH^w(t)
 \leq \mu(t) \log {t^2\over 2\mu(t)}.
 \label{equation:25} 
\end{equation}
(ii) If $\tau(t)$ is decreasing in $t\in {\cal S}$, then
\begin{equation}
 \barH^w(t)
 \leq \int_0^t x\,\tau(x)\,{\rm d}x-\mu(t)\,[1+\log \tau(t)].
 \label{equation:17} 
\end{equation}
\end{theorem}
\begin{proof}
Eq.\ (\ref{equation:25}) is an immediate consequence of (\ref{equation:24}). 
Furthermore, from (\ref{equation:7}) we have:
$$
 \barH^w(t)=-\frac{1}{F(t)}\int_0^t x\,f(x)\,\log\tau(x)\,{\rm d}x
 -\frac{1}{F(t)}\int_0^t x\,f(x)\,\log\frac{F(x)} {F(t)}\,{\rm d}x.
$$
Since $\tau(t)$ is decreasing in $t\in {\cal S}$, we have $\log\tau(x)\geq\log\tau(t)$ 
for $0<x<t$. Moreover, recalling that $\log x<x-1$ for $x>0$, we obtain:
$$
\barH^w(t)\leq -\frac{\log\tau(t)}{F(t)}\int_0^t x\,f(x)\,{\rm d}x+\frac{1}{F(t)}\int_0^t
x\,f(x)\left[\frac{F(t)}{F(x)}-1\right] {\rm d}x.
$$
By (\ref{equation:13}) and (\ref{equation:15}) we finally come to (\ref{equation:17}).
\end{proof}
\begin{remark}{\rm 
If $X$ is uniformly distributed on $(0,\nu)$, the weighted past  entropy is: 
\begin{equation}
 \barH^w(t)
 =\frac{t}{2}\log t,
 \qquad 0<t<\nu.
 \label{equation:26} 
\end{equation}
Hence, since $\mu(t)=t/2$, Eq.\ (\ref{equation:25}) is satisfied with the 
equality sign for all $t\in(0,\nu)$. 
}\end{remark}
\par
We shall now introduce two new classes of distributions based on monotonicity 
properties of the weighted entropies.
\par
\begin{definition}
A random lifetime $X$ will be said to have \\
(i) decreasing [increasing] weighted uncertainty residual life (DWURL)
[IWURL] if and only if  $H^w(t)$ is decreasing [increasing] in $t\in{\cal S}$; \\
(ii) decreasing [increasing] weighted uncertainty past life (DWUPL)
[IWUPL] if and only if  $\barH^w(t)$ is decreasing [increasing] in $t\in{\cal S}$.
\end{definition}
\par
Let $X$ be a random variable uniformly distributed on $(0,\nu)$. Since 
$H^w(t)=\frac{t+\nu}{2}\log(b-t)$, $0<t<\nu$, $X$ is DWURL if and only if 
$0<\nu\leq e$, and it can never be IWURL. Moreover, from (\ref{equation:26}) 
$X$ is DWUPL if and only if $0<\nu\leq {1\over e}$, and it can never be IWUPL. 
Finally, by virtue of (\ref{equation:27}), an exponential distribution with 
parameter $\lambda$ is DWURL (IWURL) if and only if $\lambda\geq e$ ($0<\lambda\leq e$). 
\section{Monotonic transformations}
In this section we study the weighted residual  entropy and the weighted past  
entropy under monotonic transformations. Similarly to Proposition 2.4 
of Di Crescenzo and Longobardi \cite{DiLo2002}, we have:  
\begin{theorem}\label{theorem:5}
Let $Y=\phi(X)$, with $\phi$ strictly monotonic, continuous and 
differentiable, with derivative $\phi'$. Then, for all $t\in{\cal S}$ 
\begin{equation} 
 \; H_Y^w(t)=
 \left\{
 \begin{array}{ll} 
 \!\! H^{w,\phi}\left(\phi^{-1}(t)\right)
 +\EE\left\{\phi(X)\log \phi'(X)\,|\,X>\phi^{-1}(t)\right\}, 
 & \mbox{$\phi$ strictly increasing} \\
 \hfill & \hfill \\
 \!\!
 \barH^{w,\phi}\left(\phi^{-1}(t)\right)
 +\EE\left\{\phi(X)\log[-\phi'(X)]\,|\,X\leq \phi^{-1}(t)\right\}, 
 & \mbox{$\phi$ strictly decreasing}
 \end{array}
 \right.  
 \label{equation:18} 
\end{equation}
and 
\begin{equation}
  \; \barH_Y^w(t)=
 \left\{
 \begin{array}{ll} 
 \!\! \barH^{w,\phi}\left(\phi^{-1}(t)\right)
 +\EE\left\{\phi(X)\log \phi'(X)\,|\,X\leq \phi^{-1}(t)\right\},
 & \mbox{$\phi$ strictly increasing} \\
 \hfill & \hfill \\
 \!\!
 H^{w,\phi}\left(\phi^{-1}(t)\right)
 +\EE\left\{\phi(X)\log [-\phi'(X)]\,|\,X>\phi^{-1}(t)\right\}, 
 & \mbox{$\phi$ strictly decreasing,}
 \end{array}
 \right. 
 \label{equation:19} 
\end{equation}
where
$$
 H^{w,\phi}(t)=-\frac{1}{\barF_X(t)}\int_t^{+\infty}\phi(x)\,f_X(x)\log
 \frac{f_X(x)}{\barF_X(t)}\,{\rm d}x,
$$
$$
 \barH^{w,\phi}(t)=-\frac{1}{F_X(t)}\int_t^{+\infty}\phi(x)\,f_X(x)\log
 \frac{f_X(x)}{F_X(t)}\,{\rm d}x.
$$
\end{theorem}
\begin{proof}
From (\ref{equation:6}) we have
$$
 H_Y^w(t)=-\int_t^{+\infty}y\,\frac{f_X(\phi^{-1}(y))}{P(Y>t)}\,
 \left|\frac{{\rm d}}{{\rm d}y}\phi^{-1}(y)\right|\,
 \log\left\{\frac{f_X(\phi^{-1}(y))}{P(Y>t)}\,
 \left|\frac{{\rm d}}{{\rm d}y}\phi^{-1}(y)\right|\right\}\,{\rm d}y.
$$
Let $\phi$ be strictly increasing. By setting $y=\phi(x)$, we obtain
$$
 H_Y^w(t)=-\int_{\phi^{-1}(t)}^{\phi^{-1}(+\infty)}\phi(x)\,\frac{f_X(x)}{\barF_X(\phi^{-1}(t))}
 \log\left\{\frac{f_X(x)}{\barF_X(\phi^{-1}(t))}\,
 \left|\frac{\rm d}{{\rm d}x}\phi(x)\right|^{-1}\right\}\,{\rm d}x,
$$
or
$$
 H_Y^w(t)=-\int_{\phi^{-1}(t)}^{+\infty}\frac{\phi(x)f_X(x)}{\barF_X(\phi^{-1}(t))}
 \log\frac{f_X(x)}{\barF_X(\phi^{-1}(t))}\,{\rm d}x
 +\int_{\phi^{-1}(t)}^{+\infty}\frac{\phi(x)f_X(x)}{\barF_X(\phi^{-1}(t))}
 \log\left|\frac{\rm d}{{\rm d}x}\phi(x)\right| {\rm d}x,
$$
giving the first of (\ref{equation:18}).
If $\phi$ is strictly decreasing we similarly obtain:
$$
 H_Y^w(t)=-\int_{\phi^{-1}(+\infty)}^{\phi^{-1}(t)}\frac{\phi(x)\,f_X(x)}{F_X(\phi^{-1}(t))}
 \log\frac{f_X(x)}{F_X(\phi^{-1}(t))}\,{\rm d}x
 +\int_0^{\phi^{-1}(t)}\frac{\phi(x)\,f_X(x)}{F_X(\phi^{-1}(t))}
 \log\left|\frac{\rm d}{{\rm d}x}\phi(x)\right| {\rm d}x,
$$
i.e., the second of (\ref{equation:18}). The proof of (\ref{equation:19}) is analogous.
\end{proof}
\par
According to Remark 2.3 of \cite{DiLo2002} we note that when $Y=\phi(X)$ is  
distributed as $X$ (like as for certain Pareto and beta-type distributions), 
Theorem \ref{theorem:5} yields useful identities that allow to express $H_X^w(t)$ 
and $\barH_X^w(t)$ in terms of $H^{w,\phi}\left(\phi^{-1}(t)\right)$ or 
$\barH^{w,\phi}\left(\phi^{-1}(t)\right)$, depending on the type of monotonicity of $\phi$. 
\begin{remark}{\rm
Due to Theorem \ref{theorem:5}, for all $a>0$ and $t>0$ there holds: 
$$
 H_{aX}^w(t)
 = a\,H^w\left(\frac{t}{a}\right)
 + \delta\left(\frac{t}{a}\right) a\,\log a,
$$
$$
 \barH_{aX}^w(t)
 = a\,\barH^w\left(\frac{t}{a}\right)
 + \mu\left(\frac{t}{a}\right) a\,\log a.
$$
Furthermore, for all $b>0$ and $t>b$ one has: 
$$
 H_{X+b}^w(t)
 = H^w(t-b)+b\,H(t-b), 
$$
$$
 \barH_{X+b}^w(t)
 = \barH^w(t-b)+b\,\barH(t-b). 
$$
}\end{remark}
%
\subsection*{\bf Acknowledgments}
%
This work has been performed under partial support by MIUR (PRIN 2005), 
by G.N.C.S.-INdAM and by Campania Region.
%

%
%

\begin{thebibliography}{99}
%
\bibitem{AsEb2000} 
{\sc Asadi, M. and Ebrahimi, N.} (2000), Residual entropy and its
characterizations in terms of hazard function and mean residual
life function. {\em Stat. Prob. Lett.} {\bf 49}, 263--269.
%
\bibitem{BaSk2006} 
{\sc Bartoszewicz, J. and Skolimowska, M.} (2006), Preservation of 
classes of life distributions and stochastic orders under weighting. 
{\em Stat. Prob. Lett.} {\bf 76}, 587--596.
%
\bibitem{BeGu1968} 
{\sc Belis, M. and Guia\c{s}u, S.} (1968), A quantitative-qualitative
measure of information in cybernetic systems. 
{\em IEEE Trans. Inf. Th.} {\bf IT-4}, 593--594.
%
\bibitem{BeGuNaRu2001} 
{\sc Belzunce, F., Guillamon, A., Navarro, J. and Ruiz, J.M.}
(2001), Kernel estimation of residual entropy. 
{\em Commun. Statist. -- Theory Meth.} {\bf 30}, 1243--1255. 
%
\bibitem{NaBeRuDe2004} 
{\sc Belzunce, F., Navarro, J., Ruiz, J.M. and del Aguila, Y.} (2004),
Some results on residual entropy function.
{\em Metrika} {\bf 59}, 147--161. 
%
\bibitem{DiLo2002}  
{\sc Di Crescenzo, A. and Longobardi, M.} (2002),
Entropy-based measure of uncertainty in past lifetime distributions. 
{\em J. Appl. Prob.} {\bf 39}, 434--440.
%
\bibitem{DiLo04}  
{\sc Di Crescenzo, A. and Longobardi, M.} (2004),
A measure of discrimination between past lifetime distributions. 
{\em Stat. Prob. Lett.\/} {\bf 67}, 173--182.
%
\bibitem{DiLoNa04}  
{\sc Di Crescenzo, A., Longobardi M. and Nastro, A.} (2004), 
On the mutual information between residual lifetime distributions. 
In: {\em Cybernetics and Systems 2004}, Vol.~1 (R.\ Trappl, Ed.), 
p.~142--145. Austrian Society for Cybernetic Studies, Vienna. 
ISBN 3-85206-169-5.
%
\bibitem{Eb96}  
{\sc Ebrahimi, N.} (1996),
How to measure uncertainty in the residual life time distribution.
{\em Sankhy\"a, Ser. A} {\bf 58}, 48--56.
%
\bibitem{Eb97}  
{\sc Ebrahimi, N.} (1997),
Testing whether lifetime distribution is decreasing uncertainty.
{\em J. Stat. Plann. Infer.} {\bf 64}, 9--19.
%
\bibitem{Eb2000}  
{\sc Ebrahimi, N.} (2000),
The maximum entropy method for lifetime distributions.
{\em Sankhy\"a, Ser. A} {\bf 62}, 236--243.
%
\bibitem{EbKi96a} 
{\sc Ebrahimi, N. and Kirmani, S.N.U.A.} (1996),
Some results on ordering of survival functions through uncertainty.
{\em Stat. Prob. Lett.} {\bf 29}, 167--176.
%
\bibitem{EbKi96b} 
{\sc Ebrahimi, N. and Kirmani, S.N.U.A.} (1996),
A measure of discrimination between two residual lifetime distributions 
and its applications. {\em Ann. Inst. Stat. Math.} {\bf 48}, 257--265.
%
\bibitem{EbPe95}  
{\sc Ebrahimi, N. and Pellerey, F.} (1995),
New partial ordering of survival functions based on the notion of uncertainty.
{\em J. Appl. Prob.} {\bf 32}, 202--211.
%
\bibitem{GiCr01}  
{\sc Giorgi, G.M. and Crescenzi, M.} (2001),
A look at the Bonferroni inequality measure in a reliability framework.
{\em Statistica} {\bf LXI} (4), 571--583.
%
\bibitem{Gu86}  
{\sc Guia\c{s}u, S.} (1986),
Grouping data by using the weighted entropy.
{\em J. Statist. Plann. Infer} {\bf 15}, 63--69.
%
\bibitem{JoGl04}  
{\sc Johnson, D.H. and Glantz, R.M.} (2004),
When does interval coding occur? 
{\em Neurocomputing} {\bf 59-60}, 13--18.
%
\bibitem{Lo76}  
{\sc Longo, G.} (1976), A noiseless coding theorem for sources having 
utilities. {\em SIAM  J. Appl. Math.} {\bf 30}, 739--748. 
%
\bibitem{NaDeRu2001}  
{\sc Navarro, J., del Aguila, Y. and Ruiz, J.M.} (2001),
Characterizations through reliability measures from weighted distributions.
{\em Stat. Papers} {\bf 42}, 395--402. 
%
\bibitem{Sh48}  
{\sc Shannon, C.E.} (1948),
A mathematical theory of communication.
{\em Bell System Tech. J.} {\bf 27}, 279--423.
%
\bibitem{Ta90}  
{\sc Taneja, I.J.} (1990),
On generalized entropy with applications.
In: {\em Lectures in Applied Mathematics and Informatics}
(L.M.\ Ricciardi, Ed.), p.~107--169.
Manchester Univ.\ Press, Manchester. 
%
\bibitem{Wi48}  
{\sc Wiener, N.} (1948, 2nd Ed.\ 1961),
{\em Cybernetics}, The MIT Press and Wiley, New York.
%
\end{thebibliography}
\end{document}